\documentclass[aos,MSNbibl,citesort,dvips]{arximspdf}

\doi{10.1214/12-AOS981} 
\referstodoi{10.1214/11-AOS949}
\volume{40}
\issue{4}
\pubyear{2012}
\firstpage{1978}
\lastpage{1983}

\makeatletter

\newcommand{\vertex}{V}
\newcommand{\trace}{\operatorname{trace}}
\newcommand{\Shat}{\widehat{S}}
\newcommand{\Lhat}{\widehat{L}}
\newcommand{\graph}{G}
\newcommand{\edge}{E}
\newcommand{\pdim}{p}
\newcommand{\numobs}{n}
\newcommand{\degmax}{d}
\newcommand{\real}{\mathbb{R}}
\newcommand{\matsnorm}[2]{|\!|\!| #1 | \! | \!|_{{#2}}}
\newcommand{\defn}{:=}
\newcommand{\EmpCov}{\widehat{\Sigma}^\numobs}
\makeatother

\begin{document}
\begin{frontmatter}

\title{Discussion: Latent variable graphical model selection via convex optimization}
\runtitle{Comment}

\begin{aug}
\author[A]{\fnms{Martin J.} \snm{Wainwright}\corref{}\ead[label=e1]{wainwrig@stat.berkeley.edu}}
\runauthor{M. J. Wainwright}
\affiliation{University of California at Berkeley}
\address[A]{Department of Statistics\\
University of California at Berkeley\\
421 Evans Hall\\
Berkeley, California 94720\\
USA\\
\printead{e1}} 
\end{aug}

\received{\smonth{2} \syear{2012}}



\end{frontmatter}
%

\section{Introduction}

It is my pleasure to congratulate the authors for an innovative and
inspiring piece of work. Chandrasekaran, Parrilo and Willsky
(hereafter CPW) have come up with a novel approach, combining ideas
from convex optimization and algebraic geometry, to the long-standing
problem of Gaussian graphical model selection with latent variables.
Their method is intuitive and simple to implement, based on solving a
convex log-determinant program with suitable choices of
regularization. In addition, they establish a number of attractive
theoretical guarantees that hold under high-dimensional scaling,
meaning that the graph size $\pdim$ and sample size $\numobs$ are
allowed to grow simultaneously.

\subsection{Background}
Recall that an undirected graphical model (also known as a Markov
random field) consists of a family of probability distributions that
factorize according to the structure of undirected graph \mbox{$\graph
= (\vertex, \edge)$}. In the multivariate Gaussian case, the
factorization translates into a sparsity assumption on the inverse
covariance or precision matrix \cite{Speed86}. In particular, given a
multivariate Gaussian random vector $(X_1, \ldots, X_\pdim)$ with
covariance matrix $\Sigma$, it is said to be Markov with respect to
the graph $\graph$ if its precision matrix $K = \Sigma^{-1}$ has
zeroes for each distinct pair of indices $(j,k)$ \emph{not} in the
edge set $\edge$ of the graph. Consequently, the sparsity pattern of
the inverse covariance $K$ encodes the edge structure of the graph.
The goal of Gaussian graphical model selection is to determine this
unknown edge structure, and hence the sparsity pattern of the inverse
covariance matrix. It can also be of interest to estimate the
matrices $K$ or $\Sigma$, for instance, in the Frobenius or
$\ell_2$-operator norm sense. In recent years, under the assumption
that all entries of $X$ are fully observed, a number of practical
methods have been proposed and shown to perform well under
high-dimensional scaling
(e.g., \cite{Meinshausen06,Rot09,RavWaiRasYu2011,CaiZho12}).

Chandrasekaran et al. tackle a challenging extension of this problem,
in which one observes only $\pdim$ coordinates of a larger $\pdim+
h$ dimensional Gaussian random vector. In this case, the $\pdim
\times\pdim$ precision matrix $K$ of the observed components need not
be sparse, but rather, by an application of the Schur complement
formula, can be written as the difference $K = S^* - L^*$. The first
matrix $S^*$ is sparse, whereas the second matrix $L^*$ is not sparse
(at least in general), but has rank at most $h$, corresponding to the
number of latent or hidden variables. Consequently, the problem of
latent Gaussian graphical model selection can be cast as a form of
\emph{matrix decomposition}, involving a splitting of the precision
matrix into sparse and low-rank components. Based on this nice
insight, CPW propose a natural \mbox{$M$-estimator} for this problem,
based on minimizing a regularized form of the (negative) log
likelihood for a multivariate Gaussian, where the elementwise
$\ell_1$-norm is used as a proxy for sparsity, and the nuclear or
trace norm as a proxy for rank. Overall, the method is based on the
convex program
%
\begin{eqnarray}\label{EqnMest}
(\Shat, \Lhat)  \in\arg\min \bigl\{ -\ell(S -L; \EmpCov) +
\lambda_n  \bigl( \gamma\|S\|_1 + \trace(L)  \bigr)  \bigr\}\\
\eqntext{\mbox{such that $S \succeq L \succeq0$},}
\end{eqnarray}
where $\ell(S-L; \EmpCov)$ is the Gaussian log-likelihood as a
function of the precision matrix $S - L$ and the empirical covariance
matrix $\EmpCov$ of the observed variables.

\subsection{Sharpness of rates} On one hand, the paper provides
attractive guarantees on the procedure (\ref{EqnMest})---namely, that
under suitable incoherence conditions (to be discussed below) and a
sample size $\numobs\succsim\pdim$, the method is guaranteed with
high probability: (a) to correctly recover the signed support of the
sparse matrix $S^*$, and hence the full graph structure; (b) to
correctly recover the rank of the component $L^*$, and hence the
number of latent variables; and (c) to yield operator norm consistency
of the order $\sqrt{\frac{\pdim}{\numobs}}$. The proof itself
involves a clever use of the primal-dual witness
method \cite{RavWaiRasYu2011}, in which one analyzes an $M$-estimator
by constructing a primal solution and an associated dual pair, and
uses the construction to show that the optimum has desired properties
(in this case, support and rank recovery) with high probability. A
major challenge, not present in the simpler problem without latent
variables, is dealing with the potential nonidentifiability of the
matrix decomposition problem (see below for further discussion); the
authors overcome this challenge via a delicate analysis of the tangent
spaces associated with the sparse and low-rank components.

On the other hand, the scaling $\numobs\succsim\pdim$ is quite
restrictive, at least in comparison to related results without latent
variables. To provide a concrete example, consider a Gaussian
graphical model with maximum degree $\degmax$. For any such graph,\vadjust{\goodbreak}
again under a set of so-called incoherence or irrepresentability
conditions, the neighborhood-based selection of approach of
Meinshausen and B\"{u}hlmann \cite{Meinshausen06} can be shown to
correctly specify the graph structure with high probability based on
$\numobs\succsim\degmax\log\pdim$ samples. Moreover, under a
similar set of assumptions, Ravikumar et al. \cite{RavWaiRasYu2011}
show that the $\ell_1$-regularized Gaussian MLE returns an estimate of
the precision matrix with operator norm error of the order
$\sqrt{\frac{\degmax^2 \log\pdim}{\numobs}}$. Consequently, whenever
the maximum degree $\degmax$ is significantly smaller than the
dimension, results of this type allow for the sample size $\numobs$ to
be much smaller than $\pdim$. This discrepancy---as to whether or not
the sample size can be smaller than the dimension---thus raises some
interesting directions for future work. More precisely, one wonders
whether or not the CPW analysis might be sharpened so as to reduce the
sample size requirements. Possibly this might require introducing
additional structure in the low-rank matrix. From the other direction,
an alternative approach would be to develop minimax lower bounds on
latent Gaussian model selection, for instance, by using
information-theoretic techniques that have been exploited in related
work on model/graph selection and covariance estimation
(e.g., \cite{Wainwright09info,CaiZho12,SanWai12}).

\subsection{Relaxing assumptions}
The CPW analysis also imposes lower bounds on the minimum absolute
values of the nonzero entries in $S^*$, as well as the minimum
nonzero singular values of $L^*$---both must scale as
$\Omega(\sqrt{\frac{\pdim}{\numobs}})$. Clearly, some sort of lower
bound on these quantities is necessary in order to establish exact
recovery guarantees, as in the results (a) and (b) paraphrased above.
It is less clear whether lower bounds of this order are the weakest
possible, and if not, to what extent they can be relaxed. For
instance, again in the setting of Gaussian graph selection without
latent variables \cite{Meinshausen06,RavWaiRasYu2011}, the minimum
values are typically allowed to be as small as
$\Omega(\sqrt{\frac{\log\pdim}{\numobs}})$. More broadly, in many
applications, it might be more natural to assume that the data is not
actually drawn from a sparse graphical model, but rather can be
well-approximated by such a model. In such settings, although exact
recovery guarantees would no longer be feasible, one would like to
guarantee that a given method, either the
$M$-estimator (\ref{EqnMest}) or some variant thereof, can recover all
entries of $S^*$ with absolute value above a given threshold, and/or
estimate the number of eigenvalues of $L^*$ above a (possibly
different) threshold. Such guarantees are possible for ordinary
Gaussian graph selection, where it is known that $\ell_1$-based
methods will recover all entries with absolute values above the
regularization parameter \cite{Meinshausen06,RavWaiRasYu2011}.

The CPW analysis also involves various types of incoherence conditions
on the matrix decomposition. As noted by the authors, some of these
assumptions are related to the incoherence or irrepresentability
conditions imposed in past work on ordinary Gaussian graph
selection \cite{Meinshausen06,Zhao06,RavWaiRasYu2011}; others are
unique to\vadjust{\goodbreak} the latent problem, since they are required to ensure
identifiability (see discussion below). It seems worthwhile to
explore which of these incoherence conditions are artifacts of a
particular methodology and which are intrinsic to the problem. For
instance, in the case of ordinary Gaussian graph selection, there are
problems for which the neighborhood-based Lasso~\cite{Meinshausen06}
can correctly recover the graph while the $\ell_1$-regularized
log-determinant approach \cite{RavWaiRasYu2011,Mei08} cannot. Moreover,
there are problems for which, with the same order of sample size, the
neighborhood-based Lasso will fail whereas an oracle method will
succeed \cite{Wainwright09info}. Such differences demonstrate that
certain aspects of the incoherence conditions are artifacts of
$\ell_1$-relaxations. In the context of latent Gaussian graph
selection, these same issues remain to be explored. For instance, are
there alternative polynomial-time methods that can perform latent
graph selection under milder incoherence conditions? What conditions
are required by an oracle-type approach---that is, involving exact
cardinality and rank constraints?

\subsection{Toward partial identifiability}
On the other hand, certain types of incoherence conditions are clearly
intrinsic to the problem. Even at the population level, it is clearly
not possible in general to identify the components $(S^*, L^*)$ based
on observing only the sum $K = S^* - L^*$. A major contribution of the
CPW paper, building from their own pioneering work on matrix
decompositions \cite{Chand09}, is to provide sufficient conditions on
the pair $(S^*, L^*)$ that ensure identifiability. These sufficient
conditions are based on a detailed analysis of the algebraic structure
of the spaces of sparse and low-rank matrices, respectively.

In a statistical setting, however, most models are viewed as
approximations to reality. With this mindset, it could be interesting
to consider matrix decompositions that satisfy a weaker notion of
partial identifiability. To provide a concrete illustration, suppose
that we begin with a matrix pair $(S^*, L^*)$ that is identifiable
based on observing the difference $K = S^* - L^*$. Now imagine that
we perturb $K$ by a matrix that is both sparse and low-rank---for
instance, a matrix of the form $E = z z^T$ where $z$ is a sparse
vector. If we then consider the perturbed matrix $\widetilde{K} :
= K + \delta E  =  S^* - L^* + \delta E$ for some suitably small
parameter $\delta$, the matrix decomposition is longer identifiable.
In particular, at the two extremes, we can choose between the
decompositions $\widetilde{K} = (S^* + \delta E) - L^*$, where the
matrix $(S^* + \delta E)$ is sparse, or the decomposition
$\widetilde{K} = S^* - (L^* - \delta E)$, where the matrix $L^* -
\delta E$ is low-rank. Note that this nonidentifiability holds
regardless of how small we choose the scalar $\delta$. However, from
a more practical perspective, if we relax our requirement of exact
identification, then such a perturbation need not be a concern as long
as $\delta$ is relatively small. Indeed, one might expect that it
should be possible to recover estimates of the pair $(S^*, L^*)$ that
are accurate up to an error proportional to $\delta$.

In some of our own recent work \cite{AgaNegWai11}, we have provided
such guarantees for a related class of noisy matrix decomposition\vadjust{\goodbreak}
problems. In particular, we consider the observation
model\footnote{Here we follow the notation of the CPW paper for the
sparse and low-rank components.}
%
\begin{equation}\label{EqnNoisy}
Y  = \mathfrak{X}(S^* - L^*) + W,
\end{equation}
where $\mathfrak{X}\dvtx \real^{\pdim\times\pdim} \rightarrow
\real^{\numobs_1 \times\numobs_2}$ is a known linear operator and $W
\in\real^{\numobs_1 \times\numobs_2}$ is a noise matrix. In the
simplest case, $\mathfrak{X}$ is simply the identity operator.
Observation models of this form (\ref{EqnNoisy}) arise in robust PCA,
sparse factor analysis, multivariate regression and robust covariance
estimation.

Instead of enforcing incoherence conditions sufficient for
identifiability, the analysis is performed under related but milder
conditions on the interaction between $S^*$ and $L^*$. For instance,
one way of controlling the radius of nonidentifiability is via
control on the ``spikiness'' of the low-rank component, as measured by
the ratio $\alpha(L^*) \defn\frac{\pdim
\|L^*\|_\infty}{\matsnorm{L^*}{F}}$, where $\|\cdot\|_\infty$
denotes the elementwise absolute maximum and $\matsnorm{\cdot}{F}$
denotes the Frobenius norm. For any nonzero $\pdim$-dimensional
matrix, this spikiness ratio ranges between~$1$ and~$\pdim$:
\begin{itemize}
\item On one hand, it achieves its minimum value by a matrix that has
all its entries equal to the same nonzero constant (e.g., $L^* = 1
1^T$, where $1 \in\real^\pdim$ is a vector of all ones).
\item On the other hand, the maximum is achieved by a matrix that
concentrates all its mass in a single position (e.g., $L^* = e_1
e_1^T$, where $e_1 \in\real^\pdim$ is the first canonical basis
vector).
\end{itemize}
Note that it is precisely this latter type of matrix that is
troublesome in sparse plus low-rank matrix decomposition, since it is
simultaneously sparse \emph{and} low-rank. In this way, the spikiness
ratio limits the effect of such troublesome instances, thereby
bounding the radius of nonidentifiability of the model. The
paper \cite{AgaNegWai11} analyzes an $M$-estimator, also based on
elementwise $\ell_1$ and nuclear norm regularization, for estimating
the pair $(S^*, L^*)$ from the noisy observation
model (\ref{EqnNoisy}). The resulting error bounds involve both terms
arising from the (possibly stochastic) noise matrix $W$ and
additional terms associated with the radius of nonidentifiability.

The same notion of partial identifiability is applicable to latent
Gaussian graph selection. Accordingly, it seems worthwhile to explore
whether similar techniques can be used to obtain error bounds with a
similar form---one component associated with the stochastic noise
(induced by sampling), and a second deterministic component.
Interestingly, under the scaling $\numobs\succsim\pdim$ assumed in
the CPW paper, the empirical covariance matrix $\EmpCov$ will be
invertible with high probability and, hence, it can be cast as an
observation model of the form (\ref{EqnNoisy})---namely, we can write
$(\EmpCov)^{-1} = S^* - L^* + W$, where the noise matrix $W$ is
induced by sampling.

\subsection{Extensions to non-Gaussian variables} A final more speculative
yet intriguing question is whether the techniques of CPW can be
extended to graphical models involving non-Gaussian variables, for
instance, those with binary or multinomial variables for a start. The
main complication here is that factorization and conditional
independence properties for non-Gaussian variables do not translate
directly into sparsity of the inverse covariance matrix. Nonetheless,
it might be possible to reveal aspects of this factorization by some
type of spectral analysis, in which context related matrix-theoretic
approaches could be brought to bear. Overall, we should all be
thankful to Chandrasekaran, Parillo and Willsky for their innovative
work and the exciting line of questions and possibilities that it has
raised for future research.


%

\printaddresses


\begin{thebibliography}{12}

\bibitem{AgaNegWai11}
\begin{bmisc}[auto:STB|2012/05/30|10:51:56]
\bauthor{\bsnm{Agarwal},~\bfnm{A.}\binits{A.}},
  \bauthor{\bsnm{Negahban},~\bfnm{S.}\binits{S.}} \AND
  \bauthor{\bsnm{Wainwright},~\bfnm{M.~J.}\binits{M.~J.}}
(\byear{2012}).
\bhowpublished{Noisy matrix decomposition via convex relaxation: Optimal rates
  in high dimensions. \textit{Ann. Statist.}
  \textbf{40} 1171--1197.}
\bptok{imsref}%
\end{bmisc}
\endbibitem

\bibitem{CaiZho12}
\begin{bmisc}[auto:STB|2012/05/30|10:51:56]
\bauthor{\bsnm{Cai},~\bfnm{T.}\binits{T.}} \AND
  \bauthor{\bsnm{Zhou},~\bfnm{H.}\binits{H.}}
(\byear{2012}).
\bhowpublished{Minimax estimation of large covariance matrices under $\ell
  _1$-norm. \emph{Statistica Sinica}. To appear}.
\bptok{imsref}%
\end{bmisc}
\endbibitem

\bibitem{Chand09}
\begin{barticle}[mr]
\bauthor{\bsnm{Chandrasekaran},~\bfnm{Venkat}\binits{V.}},
  \bauthor{\bsnm{Sanghavi},~\bfnm{Sujay}\binits{S.}},
  \bauthor{\bsnm{Parrilo},~\bfnm{Pablo~A.}\binits{P.~A.}} \AND
  \bauthor{\bsnm{Willsky},~\bfnm{Alan~S.}\binits{A.~S.}}
(\byear{2011}).
\btitle{Rank-sparsity incoherence for matrix decomposition}.
\bjournal{SIAM J. Optim.}
\bvolume{21}
\bpages{572--596}.
\bid{doi={10.1137/090761793}, issn={1052-6234}, mr={2817479}}
\bptnote{check year}%
\bptok{imsref}%
\end{barticle}
\endbibitem

\bibitem{Mei08}
\begin{barticle}[mr]
\bauthor{\bsnm{Meinshausen},~\bfnm{Nicolai}\binits{N.}}
(\byear{2008}).
\btitle{A note on the {L}asso for {G}aussian graphical model selection}.
\bjournal{Statist. Probab. Lett.}
\bvolume{78}
\bpages{880--884}.
\bid{doi={10.1016/j.spl.2007.09.014}, issn={0167-7152}, mr={2398362}}
\bptok{imsref}%
\end{barticle}
\endbibitem

\bibitem{Meinshausen06}
\begin{barticle}[mr]
\bauthor{\bsnm{Meinshausen},~\bfnm{Nicolai}\binits{N.}} \AND
  \bauthor{\bsnm{B{\"u}hlmann},~\bfnm{Peter}\binits{P.}}
(\byear{2006}).
\btitle{High-dimensional graphs and variable selection with the lasso}.
\bjournal{Ann. Statist.}
\bvolume{34}
\bpages{1436--1462}.
\bid{doi={10.1214/009053606000000281}, issn={0090-5364}, mr={2278363}}
\bptok{imsref}%
\end{barticle}
\endbibitem

\bibitem{RavWaiRasYu2011}
\begin{barticle}[mr]
\bauthor{\bsnm{Ravikumar},~\bfnm{Pradeep}\binits{P.}},
  \bauthor{\bsnm{Wainwright},~\bfnm{Martin~J.}\binits{M.~J.}},
  \bauthor{\bsnm{Raskutti},~\bfnm{Garvesh}\binits{G.}} \AND
  \bauthor{\bsnm{Yu},~\bfnm{Bin}\binits{B.}}
(\byear{2011}).
\btitle{High-dimensional covariance estimation by minimizing {$\ell\sb
  1$}-penalized log-determinant divergence}.
\bjournal{Electron. J. Stat.}
\bvolume{5}
\bpages{935--980}.
\bid{doi={10.1214/11-EJS631}, issn={1935-7524}, mr={2836766}}
\bptok{imsref}%
\end{barticle}
\endbibitem

\bibitem{Rot09}
\begin{barticle}[mr]
\bauthor{\bsnm{Rothman},~\bfnm{Adam~J.}\binits{A.~J.}},
  \bauthor{\bsnm{Bickel},~\bfnm{Peter~J.}\binits{P.~J.}},
  \bauthor{\bsnm{Levina},~\bfnm{Elizaveta}\binits{E.}} \AND
  \bauthor{\bsnm{Zhu},~\bfnm{Ji}\binits{J.}}
(\byear{2008}).
\btitle{Sparse permutation invariant covariance estimation}.
\bjournal{Electron. J. Stat.}
\bvolume{2}
\bpages{494--515}.
\bid{doi={10.1214/08-EJS176}, issn={1935-7524}, mr={2417391}}
\bptok{imsref}%
\end{barticle}
\endbibitem

\bibitem{SanWai12}
\begin{bmisc}[auto:STB|2012/05/30|10:51:56]
\bauthor{\bsnm{Santhanam},~\bfnm{N.~P.}\binits{N.~P.}} \AND
  \bauthor{\bsnm{Wainwright},~\bfnm{M.~J.}\binits{M.~J.}}
(\byear{2012}).
\bhowpublished{Information-theoretic limits of selecting binary graphical
  models in high dimensions. \emph{IEEE Trans. Inform. Theory} \textbf{58} 4117--4134}.
\bptok{imsref}%
\end{bmisc}
\endbibitem

\bibitem{Speed86}
\begin{barticle}[mr]
\bauthor{\bsnm{Speed},~\bfnm{T.~P.}\binits{T.~P.}} \AND
  \bauthor{\bsnm{Kiiveri},~\bfnm{H.~T.}\binits{H.~T.}}
(\byear{1986}).
\btitle{Gaussian {M}arkov distributions over finite graphs}.
\bjournal{Ann. Statist.}
\bvolume{14}
\bpages{138--150}.
\bid{doi={10.1214/aos/1176349846}, issn={0090-5364}, mr={0829559}}
\bptok{imsref}%
\end{barticle}
\endbibitem

\bibitem{Wainwright09info}
\begin{barticle}[mr]
\bauthor{\bsnm{Wainwright},~\bfnm{Martin~J.}\binits{M.~J.}}
(\byear{2009}).
\btitle{Information-theoretic limits on sparsity recovery in the
  high-dimensional and noisy setting}.
\bjournal{IEEE Trans. Inform. Theory}
\bvolume{55}
\bpages{5728--5741}.
\bid{doi={10.1109/TIT.2009.2032816}, issn={0018-9448}, mr={2597190}}
\bptok{imsref}%
\end{barticle}
\endbibitem

\bibitem{Zhao06}
\begin{barticle}[mr]
\bauthor{\bsnm{Zhao},~\bfnm{Peng}\binits{P.}} \AND
  \bauthor{\bsnm{Yu},~\bfnm{Bin}\binits{B.}}
(\byear{2006}).
\btitle{On model selection consistency of {L}asso}.
\bjournal{J. Mach. Learn. Res.}
\bvolume{7}
\bpages{2541--2563}.
\bid{issn={1532-4435}, mr={2274449}}
\bptok{imsref}%
\end{barticle}
\endbibitem

\end{thebibliography}
\end{document}